\pdfoutput=1
\RequirePackage{ifpdf}
\ifpdf % We are running pdfTeX in pdf mode
\documentclass[pdftex]{sigma}
\else
\documentclass{sigma}
\fi

\DeclareMathOperator{\End}{End}

\newcommand{\g}[1]{t_{#1}^{(-1)}}

\newtheorem{Theorem}{Theorem}[section]
\newtheorem{claim}[Theorem]{Claim}
\newtheorem{maintheorem}[Theorem]{Main Theorem}
\newtheorem{Corollary}[Theorem]{Corollary}
\newtheorem{Lemma}[Theorem]{Lemma}
\newtheorem{Proposition}[Theorem]{Proposition}
\newtheorem{problem}[Theorem]{Problem}
{\theoremstyle{definition}
\newtheorem{Definition}[Theorem]{Definition}
\newtheorem{Remark}[Theorem]{Remark}}

\begin{document}

\allowdisplaybreaks

\renewcommand{\thefootnote}{$\star$}

\renewcommand{\PaperNumber}{040}

\FirstPageHeading

\ShortArticleName{Mystic Ref\/lection Groups}

\ArticleName{Mystic Ref\/lection Groups\footnote{This paper is a~contribution to the Special Issue in honor of Anatol
Kirillov and Tetsuji Miwa.
The full collection is available at \href{http://www.emis.de/journals/SIGMA/InfiniteAnalysis2013.html}
{http://www.emis.de/journals/SIGMA/InfiniteAnalysis2013.html}}}

\Author{Yuri BAZLOV~$^\dag$ and~Arkady BERENSTEIN~$^\ddag$}

\AuthorNameForHeading{Y.~Bazlov and A.~Berenstein}

\Address{$^\dag$~School of Mathematics, University of Manchester, Oxford Road, Manchester, M13 9PL, UK}
\EmailD{\href{mailto:yuri.bazlov@manchester.ac.uk}{yuri.bazlov@manchester.ac.uk}}

\Address{$^\ddag$~Department of Mathematics, University of Oregon, Eugene, OR 97403, USA}
\EmailD{\href{mailto:arkadiy@math.uoregon.edu}{arkadiy@math.uoregon.edu}}

\ArticleDates{Received December 25, 2013, in f\/inal form March 24, 2014; Published online April 04, 2014}

\Abstract{This paper aims to systematically study {\it mystic reflection groups} that emerged independently in the
paper~[\textit{Selecta Math.~(N.S.)} \textbf{14} (2009),
  325--372] by the authors and in the paper~[\textit{Algebr. Represent. Theory} \textbf{13}
  (2010), 127--158] by Kirkman, Kuzmanovich and Zhang.
A detailed analysis of this class of groups reveals that they are in a~nontrivial correspondence with the complex
ref\/lection groups $G(m,p,n)$.
We also prove that the group algebras of corresponding groups are isomorphic and classify all such groups up to
isomorphism.}

\Keywords{complex ref\/lection; mystic ref\/lection group; thick subgroups}

\Classification{16G99; 20F55; 16S80}

\renewcommand{\thefootnote}{\arabic{footnote}}
\setcounter{footnote}{0}

\vspace{-2mm}

\section{Introduction}

\looseness=-1
Let $V$ be a~complex vector space with basis $\{x_1,\ldots,x_n\}$.
Denote by~$S(V)$ the symmetric algebra of~$V$.
It is a~fundamental fact that the algebra $S(V)^{\mathbb{S}_n}$ of symmetric polynomials is isomorphic to~$S(V)$.
More generally, the Chevalley--Shephard--Todd theorem asserts that for a~f\/inite group $G\subset \mathrm{GL}(V)$, $S(V)^G$
is isomorphic to~$S(V)$ if and only if~$G$ is generated by complex ref\/lections on~$V$.

In a~remarkable paper~\cite{KKZ}, Kirkman, Kuzmanovich and Zhang solved the following problem:

\begin{problem}\label{prob:kkz}
Given a~complex matrix ${\mathbf q}=(q_{ij},1\le i,j\le n)$ with $q_{ij}q_{ji}=1$, $q_{ii}=1$, let $S_{\mathbf q}(V)$ be
the algebra generated by $V$ subject to the relations $x_ix_j=q_{ij}x_jx_i$ for $1\le i,j\le n$.
Find all finite groups~$G$ such that
\begin{itemize}\itemsep=0pt
\item[$(a)$] $G$ acts on the algebra $S_{\mathbf q}(V)$ by degree-preserving automorphisms;
\item[$(b)$] The fixed point
algebra $S_{\mathbf q}(V)^G$ is isomorphic to $S_{{\mathbf q}'}(V)$ for some ${\mathbf q'}$.
\end{itemize}
\end{problem}

We will refer to a~group~$G$ satisfying (a) and (b) above as a~\textit{mystic reflection group}.

Independently, in~\cite{BBselecta} we solved the following problem:

\begin{problem}
\label{prob:braided cherednik}
Classify all algebras ${\mathcal A}$ such that
\begin{itemize}\itemsep=0pt
\item
${\mathcal A}$ is generated by the space $V$, its dual $V^*$ and a~finite group $G\subset
\mathrm{GL}(V)$;
\item
${\mathcal A}$ admits a~triangular decomposition ${\mathcal A}=S_{\mathbf q}(V)\otimes \mathbb C G\otimes S_{\mathbf q}
(V^*)$ for some ${\mathbf q}$ as above;
\item
$S_{\mathbf q}(V)$ and $S_{\mathbf q}(V^*)$ are invariant under conjugation by elements of~$G$, and the restriction of
these conjugation representations to $V$ and $V^*$ is the natural action of~$G$ on these spaces;
\item
$y_jx_i-q_{ij}x_iy_j\in \mathbb C G$ for all $i$, $j$, where $\{y_1,\ldots,y_n\}$ is the basis of $V^*$ dual to
$\{x_1,\ldots,x_n\}$.
\end{itemize}
\end{problem}

Comparing~\cite[Theorem 1.1 and Corollary 5.6]{KKZ} and~\cite[Theorem 0.7]{BBselecta}, one obtains the following
surprising result:

\begin{Theorem}
A group~$G$ solves Problem~{\rm \ref{prob:braided cherednik}} if and only if it is a~mystic reflection group.
\end{Theorem}

The goal of this paper is to ``demystify'' the mystic ref\/lection groups, develop their structural theory, show that
their group algebras are isomorphic to those of complex ref\/lection groups, and to deduce Problem~\ref{prob:kkz} directly
from the classical Chevalley--Shepard--Todd theorem.

\section{Main results}
\label{sect:main}

We start with a~new notion of ``mystical equivalence'' of group actions, which is crucial for what follows.

\begin{Definition}
Let $ \mathop{\triangleright}\nolimits\colon G \times \mathcal V\to \mathcal V$,
$\mathop{\triangleright}\nolimits\nolimits'\colon G' \times \mathcal V\to \mathcal V $ be faithful actions of f\/inite
groups~$G$, respectively $G'$, on a~complex vector space~$\mathcal V$.
We say that the actions $\mathop{\triangleright}\nolimits$ and $\mathop{\triangleright}\nolimits'$ are
\textit{mystically equivalent}, if
\begin{gather*}
\rho(e_G)=\rho'(e_{G'}),
\end{gather*}
where $\rho\colon \mathbb C G \to \End_\mathbb C \mathcal V$ and $\rho'\colon \mathbb C G'\to \End_\mathbb C\mathcal V$
are the algebra homomorphisms def\/ined by the actions, and $e_G$ denotes the element $\sum\limits_{g\in G}g$ of $\mathbb C G$.
\end{Definition}

Mystical equivalence of the actions of~$G$ and $G'$ is a~strengthening of the condition that the respective spaces
$\mathcal V^G$ and $\mathcal V^{G'}$ of invariants are equal, due to the following obvious result.
\begin{Lemma}%\label{lemma:projection}
If a~$G$-action and a~$G'$-action on $\mathcal V$ are mystically equivalent, then $\mathcal V^G=\mathcal V^{G'}$.
\end{Lemma}

We will use the lemma in the situation where $\mathcal V = S(V)$ where $V$ is a~vector space over $\mathbb C$ with
a~chosen basis $\{x_1,\ldots,x_n\}$.
Throughout the paper, $n\ge 2$.
Denote by $\mathbb G_n$ the group of monomial matrices on $V$, that is, matrices in $\mathrm{GL}_n(\mathbb C)$ with
exactly $n$ non-zero entries.
In other words,
\begin{gather*}
\mathbb G_n=(\mathbb C^\times)^n\rtimes \mathbb S_n.
\end{gather*}
Here $(\mathbb C^\times)^n$ is naturally identif\/ied with the group of diagonal matrices in $\mathrm{GL}_n(\mathbb C)$
and acts on~$V$ by scaling the basis $\{x_1,\ldots,x_n\}$.
The symmetric group $\mathbb S_n$ is identif\/ied with the group of permutation matrices and acts on~$V$ by permuting the
same basis.
Note that $\mathbb G_n$ is the normalizer of the torus $(\mathbb C^\times)^n$ in $\mathrm{GL}_n(\mathbb C)$.
In particular, $\mathbb S_n$ acts on $(\mathbb C^\times)^n$ by conjugation.
When we write $tw\in \mathbb G_n$, we will imply that $t\in (\mathbb C^\times)^n$ and $w\in \mathbb S_n$; every element
of $\mathbb G_n$ can be uniquely written in this way.

Clearly, $\mathbb G_n$ is generated by $s_1,\ldots,s_{n-1}$ and $t_j^{(\zeta)}$, $1\le j\le n$, $\zeta\in \mathbb
C^\times$, where
\begin{itemize}\itemsep=0pt
\item
$s_i\in \mathbb S_n$ is the permutation of $\{x_1,\ldots,x_n\}$ which swaps $x_i$ and $x_{i+1}$;
\item
$t_j^{(\zeta)}\in (\mathbb C^\times)^n$ maps $x_k$ to $\zeta^{\delta_{jk}} x_k$.
\end{itemize}
We f\/ind it very convenient to use the linear character
\begin{gather*}
\det \colon \ \mathbb G_n \to \mathbb C^\times
\end{gather*}
of $\mathbb G_n$, which is just the restriction of the determinant character of $\mathrm{GL}_n(\mathbb C)$ to $\mathbb
G_n$.
In particular,
\begin{gather*}
\det s_i = -1,
\qquad
\det w\in\{\pm1\}\text{ for }w\in \mathbb S_n,
\qquad
\det\left(t_1^{(\zeta_1)} t_2^{(\zeta_2)} \cdots t_n^{(\zeta_n)}\right)=\zeta_1\zeta_2\cdots \zeta_n.
\end{gather*}

{\sloppy Next, we introduce \textit{two} dif\/ferent faithful actions of $\mathbb G_n$ on $S(V)$ using the natural basis
$\big\{x_1^{k_1} \cdots x_n^{k_n} :  k_1,\ldots,k_n\in \mathbb Z_{\ge0}\big\}$ of~$S(V)$.
(At the moment, we are not using any multiplication on $S(V)$.)

}

\begin{Proposition}\label{prop:automorphism}\quad
\begin{enumerate}\itemsep=0pt
\item[$(a)$] There exist $($unique$)$ faithful actions $\mathop{\triangleright}\nolimits_+$, $\mathop{\triangleright}\nolimits_-$ of
$\mathbb G_n$ on $S(V)$ such that
\begin{gather*}
s_i \mathop{\triangleright}\nolimits_\pm x_1^{k_1}\cdots x_n^{k_n} = (\pm 1)^{k_ik_{i+1}} x_1^{k_1}\cdots
x_i^{k_{i+1}}x_{i+1}^{k_i}\cdots x_n^{k_n},
\\
t_j^{(\zeta)}\mathop{\triangleright}\nolimits_\pm x_1^{k_1}\cdots x_n^{k_n} = \zeta^{k_j} x_1^{k_1}\cdots
x_n^{k_n}
\end{gather*}
for any $i=1,\ldots,n-1$, $j=1,\ldots,n$, $\zeta\in \mathbb C^\times$.
Both actions extend the defining action of~$\mathbb G_n$ on~$V$.

\item[$(b)$] The action $\mathop{\triangleright}\nolimits_+$ of $\mathbb G_n$ on $S(V)$ is compatible with the natural
commutative multiplication on~$S(V)$, in the sense that~$\mathbb G_n$ acts by automorphisms of the algebra~$S(V)$.

\item[$(c)$] The action $\mathop{\triangleright}\nolimits_-$ of $\mathbb G_n$ on $S(V)$ is compatible    with the
algebra structure $S_{\mathbf{-1}}(V)$ on $S(V)$ $($which is $S_{\mathbf q}(V)$ with $q_{ij}=-1$ for all $i\ne j)$;

\item[$(d)$] $\rho_+(\mathbb C \mathbb G_n)= \rho_-(\mathbb C \mathbb G_n)$, where $\rho_\pm\colon \mathbb C\mathbb G_n \to
\End_\mathbb C S(V)$ are the algebra homomorphisms arising from the actions $\mathop{\triangleright}\nolimits_\pm$.
Moreover, $\rho_-=\rho_+ \circ J$ for some algebra automorphism~$J$ of~$\mathbb C \mathbb G_n$.
\end{enumerate}
\end{Proposition}

\begin{Remark}
\label{rem:injectivity}
By a~powerful result on group actions on integral domains, see Corollary~\ref{cor:premet-dedekind} from the
Appendix taken with $R=\mathbb C$, $A=S(V)$, the action $\mathop{\triangleright}\nolimits_+$ is faithful, and
moreover the corresponding algebra homomorphism $\rho_+\colon \mathbb C \mathbb G_n \to \End_\mathbb C S(V)$ is
injective.
It follows from (d) that $\rho_-\colon \mathbb C \mathbb G_n\to \End_\mathbb C S_{\mathbf{-1}}(V)$ is also injective.
This fact does not readily follow from classical results.
\end{Remark}

At this point, we restrict our attention to f\/inite subgroups~$G$ of $\mathbb G_n$~-- namely, to Shephard--Todd's
imprimitive complex ref\/lection groups $G(m,p,n)$ and the groups $W_{{\mathcal C},{\mathcal C}'}$, introduced
independently in~\cite{BBselecta} and~\cite{KKZ} and def\/ined as follows.
Let $n\ge 1$ and ${\mathcal C}'\subseteq{\mathcal C}$ be two f\/inite subgroups of $\mathbb C^\times$ of orders $\frac
mp$, $m$ respectively.
Then
\begin{gather*}
G(m,p,n)=\{tw\in {\mathcal C}^n\rtimes \mathbb S_n: \det t\in{\mathcal C}'\},
\\
W_{{\mathcal C},{\mathcal C}'} = \{ tw\in {\mathcal C}^n\rtimes \mathbb S_n: \det(tw)\in {\mathcal C}'\}.
\end{gather*}
The similarity of the two def\/initions manifests itself in our f\/irst main result  where a~correspondence, $\mu$,
between the two classes of subgroups of $\mathbb G_n$ is established.

\begin{maintheorem}
\label{thm:main}
Given $n\ge 1$, an even $m\ge 2$ and a~divisor $p$ of $m$, let $G=G(m,p,n)$.
Then there exists a~unique finite subgroup $ \mu(G) \subset \mathbb G_n $ such that the restriction of
$\mathop{\triangleright}\nolimits_-$ onto $\mu(G)$ is mystically equivalent to $\mathop{\triangleright}\nolimits_+$ on~$G$.
In fact,
\begin{gather*}
\mu(G)=W_{\mathcal C,\mathcal C'},
\end{gather*}
where $|{\mathcal C}|=m$, $|{\mathcal C}'|=\frac mp$.
\end{maintheorem}

The def\/inition of the group $\mu(G)$ suggests that the invariants of $\mu(G(m,p,n))=W_{{\mathcal C},{\mathcal C}'}$
should be viewed in the noncommutative algebra $S_{\mathbf{-1}}(V)$, where this group acts via
$\mathop{\triangleright}\nolimits_-$.
To describe these invariants, introduce the following elements in the space $S(V)$:
\begin{gather*}
p_k^{(m)} = \sum\limits_{i=1}^n x_i^{km},
\qquad
k=1,\ldots,n-1,
\qquad
r^{(l)} = x_1^{l}x_2^{l}\cdots x_n^{l}.
\end{gather*}
The classical result of Shephard--Todd and Chevalley asserts that the subalgebra $S(V)^{G(m,p,n)}$ of~$S(V)$ (with
respect to the natural commutative product on~$S(V)$) equals
the polynomial algebra $\mathbb C\big[p_1^{(m)},\ldots,p_{n-1}^{(m)}, r^{(\frac{m}{p})}\big]$.

Our mystic equivalence construction immediately leads to the following result, f\/irst obtained by Kirkman, Kuzmanovich
and Zhang as a~key ingredient in the classif\/ication theorem~\cite[Theorem~1.1]{KKZ}.

\begin{Theorem}
In the notation of Theorem~{\rm \ref{thm:main}}, let $m$ be even.
Then in the algebra $S_{\mathbf{-1}}(V)$,
\begin{enumerate}\itemsep=0pt
\item[$(a)$] the elements $p_1^{(|\mathcal C|)},\ldots,p_{n-1}^{(|\mathcal C|)}, r^{(|\mathcal C'|)}$ are pairwise commuting
invariants of the group $W_{\mathcal C,\mathcal C'}$;

\item[$(b)$] $S_{\mathbf{-1}}(V)^{W_{\mathcal C,\mathcal C'}} = \mathbb C\big[p_1^{(|\mathcal C|)},\ldots,p_{n-1}^{(|\mathcal C|)},r^{(|\mathcal C'|)}\big]$.
\end{enumerate}
\end{Theorem}

\begin{Remark}%\label{rem:EKL}
This result shows that the condition that $m$ is even in Theorem~\ref{thm:main} is important: the symmetric group
$\mathbb S_n=G(1,1,n)$ does not have a~mystical counterpart, and indeed the correspondence $\mu$ cannot be extended to
groups $G(m,p,n)$ where $m$ is odd.
This happens because the space of the invariants of $G(m,p,n)$ in $S(V)$ is not closed under the multiplication in
$S_{\mathbf{-1}}(V)$, and therefore cannot be the space of invariants of a~group acting by automorphisms of
$S_{\mathbf{-1}}(V)$.
\end{Remark}

The groups $G=G(m,p,n)$ and $\mu(G)$ are of the same order $\frac{m^n n!}{p}$, and, informally, they ``look very similar''.
Our next main result makes this informal statement more precise.
We keep the notation from Theorem~\ref{thm:main} and denote by $R$ the ring $\mathbb
Z\big[\frac{1+\mathbf{i}}{2}\big]\subset\mathbb C$.
\begin{maintheorem}
\label{thm:main2}
For all~$G$ as in Theorem~{\rm \ref{thm:main}}, the group rings $RG$ and $R\mu(G)$ are isomorphic.
In particular, $ \mathbb C G \cong \mathbb C \mu(G) $.
\end{maintheorem}

This theorem is rather nontrivial because the groups~$G$ and $\mu(G)$ are often not isomorphic as abstract groups.
Indeed, it was shown in~\cite[Example 7.3]{KKZ} that the group $G=G(2,2,n)$ is not isomorphic to its mystic counterpart
$\mu(G)=W_{\{\pm1\},\{1\}}$ for all even $n$.
We generalize this observation and give a~complete list of cases where~$G$ is not isomorphic to $\mu(G)$.

\begin{Theorem}
\label{thm:not iso}
In the notation of Theorem~{\rm \ref{thm:main}}, let $G=G(m,p,n)$ with $m$ even.
Then the groups~$G$ and $\mu(G)$ are not
isomorphic as abstract groups, if and only if $n$ is even and $\frac mp$ is odd.
\end{Theorem}

We go further than this and classify all groups of the form $G=G(m,p,n)$ and $\mu(G)$ up to isomorphism.
We need the following useful notion.

\begin{Definition}
We say that a~subgroup~$G$ of a~semidirect product $T\rtimes H$ is \textit{thick} if
\begin{itemize}
\item
$\pi(G)=H$, where $\pi\colon T\rtimes H \to H$ is the canonical projection onto the second factor;
\item
$G$ is normal in in $T\rtimes H$.
\end{itemize}
\end{Definition}
It is not dif\/f\/icult to see that all the groups from Theorem~\ref{thm:main} are thick subgroups of $G(m,1,n) ={\mathcal
C}^n\rtimes \mathbb S_n$.
It turns out that a~converse is also true.

\begin{Theorem}
\label{thm:thick}
Let ${\mathcal C}$ be the subgroup of $\mathbb C^\times$ of order $m$.
Then every thick subgroup of ${\mathcal C}^n \rtimes \mathbb S_n$ is of the form $G(m,p,n)$ or $($if $m$ is even$)$
$W_{{\mathcal C},{\mathcal C}'}$, and in particular, is a~mystic reflection group.
\end{Theorem}

The following completes the classif\/ication of thick subgroups of all $G(m,1,n)$ up to isomorphism.
We keep the notation used in the preceding theorems.

\begin{Theorem}
\label{thm:classification}
Let $G\subseteq G(m,1,n)$, $G'\subseteq G(m',1,n')$ be thick subgroups.
\begin{enumerate}\itemsep=0pt
\item[$(a)$] Suppose that $n=n'$ and $G\ne G'$ in $\mathbb G_n$.
Then $G\cong G'$ if and only if $n$ is odd, $m=m'$ is even, and $\{G,G'\}=\{G(m,p,n),\mu(G(m,p,n))\}$ for some $p$ such
that $\frac{m}{p}$ is odd.

\item[$(b)$] Suppose that $n\!<\!n'$.
Then $G\cong G'$ if and only if $n\!=\!3$, $n'\!=\!4$, $G\!\in\!\{G(2,2,3),\mu(G(2,2,3))\}$, $G'=G(1,1,4)=\mathbb S_4$.
\end{enumerate}
\end{Theorem}

\begin{Remark}
It is not dif\/f\/icult to see that if $\frac{m}{p}$ is even (in the notation of Theorem~\ref{thm:main}), then $G=\mu(G)$ in
$\mathbb G_n$.
This theorem together with Theorem~\ref{thm:main2} implies that the converse is also true.
\end{Remark}

The above classif\/ication suggests the following general problem which we do not address in the present paper:

\begin{problem}
Given a~semidirect product group $T\rtimes H$ where $T$ and $H$ are finite, classify all finite thick subgroups of
$T\rtimes H$ up to isomorphism.
\end{problem}

\section{Proofs of results from Section~\ref{sect:main}}

We will repeatedly use the following straightforward technical fact about the root system of type~$A_{n-1}$, the proof
of which is left to the reader as an exercise.
\begin{claim}
\label{claim:cocycle}
Let $A$ be a~multiplicatively written abelian group.
Let $\phi_{ij}\colon \mathbb Z^n\to A$, $1\le i,j\le n$, $i\ne j$, be a~system of maps satisfying
\begin{gather*}
\phi_{ij}(w(\mathbf k)) = \phi_{w^{-1}(i),w^{-1}(j)}(\mathbf k),
\qquad
\phi_{ji}(\mathbf k)\phi_{ij}(\mathbf k)=1
\end{gather*}
for all $\mathbf k=(k_1,\ldots,k_n)\in\mathbb Z^n$, $w\in\mathbb S_n$, where $w(\mathbf k)$ stands for
$(k_{w^{-1}(1)},\ldots,k_{w^{-1}(n)})$.
Denote
\begin{gather*}
\phi_w(\mathbf k) = \prod\limits_{i<j,w(i)>w(j)} \phi_{ij}(\mathbf k).
\end{gather*}
Then:
\begin{enumerate}\itemsep=0pt
\item[$(a)$] For all $w',w\in \mathbb S_n$,
\begin{gather*}
\phi_{w'w}(\mathbf k) = \phi_{w'}(w(\mathbf k))\phi_w(\mathbf k).
\end{gather*}

\item[$(b)$] If for all $i\ne j$ and for all $\mathbf k,\mathbf{k'}\in \mathbb Z^n$ one has
\begin{gather*}
\phi_{ij}(\mathbf k+\mathbf{k'}) = a_{ij}^{k_ik_j'- k_i'k_j} \phi_{ij}(\mathbf k)\phi_{ij}(\mathbf{k'}),
\end{gather*}
then for all $w\in \mathbb S_n$ one has
\begin{gather*}
\langle \mathbf k,\mathbf{k'}\rangle \phi_{w}(\mathbf k+\mathbf{k'}) = \langle w(\mathbf k),w(\mathbf{k'})\rangle
\phi_{w}(\mathbf k)\phi_{w}(\mathbf{k'}),
\end{gather*}
where $\langle \mathbf k,\mathbf{k'}\rangle= \prod\limits_{i<j}a_{ij}^{k_i'k_j}$.
Here $a_{ij}$ are elements of $A$ such that $a_{ji}=a_{ij}$.
\end{enumerate}
\end{claim}

\subsection*{Proof of Proposition~\ref{prop:automorphism}}

Observe that the group $\mathbb G_n\subset \mathrm{GL}_n(\mathbb C)$ is generated by $\mathbb S_n$ and $(\mathbb
C^\times)^n$ subject to the semidirect product relations
\begin{gather*}
wt_j^{(\zeta)}=t_{w(j)}^{(\zeta)} w,
\qquad
\text{for all} \ \ w\in \mathbb S_n,
\quad
j\in\{1,\ldots,n\},
\quad
\zeta\in\mathbb C^\times.
\end{gather*}
In what follows, we use the abbreviation $x^{\mathbf k}$ to denote $x_1^{k_1}\cdots x_n^{k_n}$, where $\mathbf
k=(k_1,\ldots,k_n)\in \mathbb Z_{\ge 0}^n$.
For $c\in \mathbb C^\times$, def\/ine $\phi_{ij}^{(c)}\colon \mathbb Z^n\to \mathbb C^\times$ by
\begin{gather*}
\phi_{ij}^{(c)}(\mathbf k) = (-1)^{k_ik_j} c^{\overline k_i-\overline k_j},
\qquad
\text{where}
\qquad
\overline k=(k\ \mathrm{mod}\ 2)\in\{0,1\}.
\end{gather*}
Additionally, def\/ine $\phi_{ij}^{(0)}(\mathbf k)=1$ for all $\mathbf k$.
Because $\phi_{ij}^{(c^{-1})}(\mathbf k)=\phi_{ij}^{(c)}(\mathbf k)^{-1}=\phi_{ji}^{(c)}(\mathbf k)$, the system~$\phi_{ij}^{(c)}$ satisf\/ies the condition in Claim~\ref{claim:cocycle} and hence gives rise to functions
$\phi_w^{(c)}$.
The following lemma is then immediate from Claim~\ref{claim:cocycle}(a).

\begin{Lemma}
\label{lem:action}
For each $c\in\mathbb C$, the formula
\begin{gather*}
w t_1^{(\zeta_1)}\cdots t_n^{(\zeta_n)} \mathop{\triangleright}\nolimits_c x^{\mathbf k} = \phi^{(c)}_w(\mathbf k)
\left(\prod\limits_{j=1}^n \zeta_j^{k_j}\right) x^{w(\mathbf k)}
\end{gather*}
defines an action $\mathop{\triangleright}\nolimits_c$ of the group $\mathbb G_n$ on the space $S(V)$.
\end{Lemma}

The actions $\mathop{\triangleright}\nolimits_0$, $\mathop{\triangleright}\nolimits_1$ given by Lemma~\ref{lem:action}
coincide on the generators $s_i$, $t_j^{(\zeta)}$ of $\mathbb G_n$ with the actions
$\mathop{\triangleright}\nolimits_+$, $\mathop{\triangleright}\nolimits_-$ def\/ined in part~(a) of
Proposition~\ref{prop:automorphism}.
Hence part~(a) of the proposition is proved.

Denote by $\bullet_+$, respectively $\bullet_-$, the multiplication on the algebra $S(V)$, respectively
$S_{\mathbf{-1}}(V)$.
One has the following multiplication rule for monomials:
\begin{gather*}
x^{\mathbf k}\bullet_\pm x^{\mathbf{k'}} = \langle \mathbf k,\mathbf k'\rangle_\pm x^{\mathbf{k}+\mathbf{k'}}
\qquad
\text{for all}\quad\mathbf{k}, \mathbf{k'}\in \mathbb Z_{\ge 0}^n,
\end{gather*}
where $\langle \mathbf k,\mathbf k'\rangle_\pm$ is as given in Claim~\ref{claim:cocycle}(b) with all $a_{ij}=\pm1$, $i\ne j$.

It is enough to check that $w\mathop{\triangleright}\nolimits_\pm$ and
$t_j^{(\zeta)}\mathop{\triangleright}\nolimits_\pm$ are automorphisms of the respective algebra structures on $S(V)$.
We apply these actions to both sides of the multiplication rule for monomials and check that the results are equal.
This is trivial for $t_j^{(\zeta)}\mathop{\triangleright}\nolimits_\pm$.
For $w\mathop{\triangleright}\nolimits_\pm$ where $w\in \mathbb S_n$,   the equality is guaranteed by
Claim~\ref{claim:cocycle}(b) and Lemma~\ref{lem:action}, applied to functions $\phi_{ij}^+=\phi_{ij}^{(0)}$,
respectively $\phi_{ij}^-=\phi_{ij}^{(1)}$.
This proves parts (b), (c) of Proposition~\ref{prop:automorphism}.

Now let us prove part (d) of Proposition~\ref{prop:automorphism}.
Clearly, the natural $\mathbb S_n$-action on the group $(\mathbb C^\times)^ n\subset \mathbb G_n$ extends to that on the
group algebra $\mathbb C (\mathbb C^\times)^n$.

For each $c\in \mathbb C\setminus\{0\}$, $w\in \mathbb S_n$, def\/ine the element $Q_w^{(c)}\in \mathbb C (\mathbb C^\times)^n$ by
\begin{gather*}
Q_w^{(c)} = \prod\limits_{i<j:w(i)>w(j)} Q_{ij}^{(c)},
\end{gather*}
where $ Q_{ij}^{(c)} = \frac14 \big(\big(c+c^{-1}\big)\big(1-\g i\g j\big)+ \big(c-c^{-1}+2\big)\g i+\big(c^{-1}-c+2\big)\g j\big)$, $1\le i,j\le n$.
\begin{Lemma}\label{lem:Q}
For all $w',w\in \mathbb S_n$, $c\in\mathbb C^\times$, $Q_w^{(c)}Q_w^{(c^{-1})}=1$ and
$Q_{w'w}^{(c)}=w^{-1}(Q_{w'}^{(c)})\cdot Q_w^{(c)}$.
\end{Lemma}
\begin{proof}
Denote by $\mathcal F$ the algebra of all functions from $\mathbb Z^n$ to $\mathbb C$ with pointwise addition and
multiplication.
Clearly, the assignment $t_1^{(\zeta_1)}\ldots t_n^{(\zeta_n)} \mapsto (\mathbf k\mapsto \zeta_1^{k_1}\ldots
\zeta_n^{k_n})$ def\/ines a~group homomorphism $(\mathbb C^\times)^n\to \mathcal F^\times $ which extends to an injective
map $\Psi\colon \mathbb C(\mathbb C^\times)^n\hookrightarrow \mathcal F$.

It is easy to check that $\Psi(Q_{ij}^{(c)})=\phi_{ij}^{(c)}$ where $\phi_{ij}^{(c)}$ is as in Lemma~\ref{lem:action}.
In particular, $Q_{ij}^{(c)}Q_{ij}^{(c^{-1})}=1$ since $\phi_{ij}^{(c)}\phi_{ij}^{(c^{-1})}=1$ and $\Psi$ is injective.
This proves the f\/irst assertion of the lemma.
To prove the second assertion, apply $\Psi^{-1}$ to Claim~\ref{claim:cocycle}(a) and use the fact that the function
$\mathbf k\mapsto \phi(w(\mathbf k))$ is mapped by $\Psi^{-1}$ to $w^{-1}(\Psi^{-1}(\phi))$ for all functions $\phi$
from the subgroup of $\mathcal F^\times$ generated by $\{\phi_{ij}^{(c)}:i\ne j\}$.
\end{proof}

Now for each $c\in\mathbb C^\times$ def\/ine the $\mathbb C$-linear map $J_c\colon \mathbb C\mathbb G_n\to \mathbb
C\mathbb G_n$ by the formula
\begin{gather*}
J_c(wt)=wtQ^{(c)}_w,
\qquad
t\in \mathbb C(\mathbb C^\times)^n,
\qquad
w\in \mathbb S_n.
\end{gather*}
\begin{Lemma}\label{lem:J_c}
For each $n\ge 1$,
\begin{enumerate}\itemsep=0pt
\item[$(a)$] $J_c$ is an algebra automorphism of $\mathbb C \mathbb G_n$ with inverse $J_{c^{-1}}$.

\item[$(b)$] $\rho_+\circ J_c=\rho_c$, where $\rho_c:\mathbb C \mathbb G_n\to \End_\mathbb C S(V)$ is the algebra homomorphism
corresponding to $\mathop{\triangleright}\nolimits_c$.
\end{enumerate}
\end{Lemma}

\begin{proof}
On the one hand, $J_c(w't'w t)=J_c(w'w\cdot w^{-1}(t')t)=w'w \cdot w^{-1}(t')tQ_{w'w}^{(c)}$.
On the other hand,
\begin{gather*}
J_c(w't')J_c(w t) = w't'Q_{w'}^{(c)} wt Q^{(c)}_w=w'w \cdot w^{-1}(t')tw^{-1}(Q_{w'}^{(c)}) Q_w^{(c)}=J_c(w't'wt)
\end{gather*}
by the second assertion of Lemma~\ref{lem:Q}.
Hence $J_c$ is a~homomorphism of algebras.
Now{\samepage
\begin{gather*}
J_c(J_{c^{-1}}(wt))=J_c\big(wtQ_w^{(c^{-1})}\big)=wtQ_w^{(c^{-1})}Q_w^{(c)}=wt
\end{gather*}
by the f\/irst assertion of Lemma~\ref{lem:Q}.
This proves part (a) of the lemma.}

Prove (b).
In view of Lemma~\ref{lem:action}, it suf\/f\/ices to show that $Q_w^{(c)}\mathop{\triangleright}\nolimits_+ x^{\mathbf
k}=\phi_w^{(c)}(\mathbf k) x^{\mathbf k}$.
Indeed,
\begin{gather*}
Q_w^{(c)}\mathop{\triangleright}\nolimits_+ x^{\mathbf k}
= \frac14 \big(\big(c+c^{-1}\big)\big(1-(-1)^{k_i+k_j}\big)+\big(c-c^{-1}+2\big)(-1)^{k_i}+\big(c^{-1}-c+2\big)(-1)^{k_j}\big)x^{\mathbf k}
\\
\phantom{Q_w^{(c)}\mathop{\triangleright}\nolimits_+ x^{\mathbf k}}
=\phi_{ij}^{(c)}(\mathbf k)x^{\mathbf k}.
\end{gather*}

Finally,
\begin{gather*}
 J_c(wt)\mathop{\triangleright}\nolimits_+ x^{\mathbf k}
=\big(wtQ_w^{(c)}\big)\mathop{\triangleright}\nolimits_+x^{\mathbf k}
=(wt)\mathop{\triangleright}\nolimits_+\big(Q_w^{(c)}\mathop{\triangleright}\nolimits_+x^{\mathbf k}\big)
=\phi_w^{(c)}(\mathbf k)(wt)\mathop{\triangleright}\nolimits_+x^{\mathbf k} =wt\mathop{\triangleright}\nolimits_c x^{\mathbf k}. \!\!\!\tag*{\qed}
\end{gather*}
\renewcommand{\qed}{}
\end{proof}

Taking $c=1$ in Lemma~\ref{lem:J_c}(b), we settle Proposition~\ref{prop:automorphism}(d).
Proposition~\ref{prop:automorphism} is proved.

\subsection*{Proof of Theorem~\ref{thm:main}} We retain the notation from the proof of
Proposition~\ref{prop:automorphism}.
Let $G=G(m,p,n)$ as in the theorem, and denote $T:=G\cap (\mathbb C^\times)^n$.
Let $e_T=\sum\limits_{t\in T}t\in \mathbb C T\subset \mathbb C G$.
Clearly, $\g ie_T=\g 1e_T$ for all $i=1,\ldots,n$.
Hence $Q_{ij}^{(c)}e_T=\g1e_T$ for all $c\in\mathbb C^\times$, so that $Q_w^{(c)}e_T=t_1^{(\det w)}e_T$ for all $w\in
\mathbb S_n$.
Since, as sets, $G=\{wt \,|\, w\in \mathbb S_n, t\in T\}$, one has $e_G=\sum\limits_{w\in\mathbb S_n}w\cdot e_T$ and
\begin{gather*}
J_c(e_G)= \sum\limits_{w\in \mathbb S_n} J_c(w)e_T = \sum\limits_{w\in \mathbb S_n} wQ_w^{(c)}e_T =\sum\limits_{w\in
\mathbb S_n} wt_1^{(\det w)}e_T = e_{\mu(G)},
\end{gather*}
because, as sets, $\mu(G)=\big\{wt_1^{(\det w)} t \,|\, w\in \mathbb S_n, t\in T\big\}$.
Since a~subgroup $G'$ of $\mathbb G_n$ is uniquely determined by $e_{G'}\in\mathbb C \mathbb G_n$, the group $\mu(G)$ is
a~unique subgroup $G'$ of $\mathbb G_n$ such that $J_c(e_G)=e_{G'}$.

Finally, by Lemma~\ref{lem:J_c}(b), $\rho_c(e_G)=\rho_+(e_{\mu(G)})$.
Setting $c=1$ and using injectivity of $\rho_+$, see Remark~\ref{rem:injectivity}, completes the proof of
Theorem~\ref{thm:main}.

\subsection*{Proof of Theorem~\ref{thm:main2}} Let $G=G(m,p,n)$.
Consider the restriction of $J_\mathbf{i}$, $\mathbf{i}=\sqrt{-1}$ to $R G=\mathbb S_n\cdot R T$ where $T=G\cap (\mathbb
C^\times)^n=\mu(G)\cap (\mathbb C^\times)^n$.
Observe that $J_\mathbf{i}(T)=T$ and
\begin{gather*}
J_\mathbf{i}(s_i)=s_iQ_{s_i}^{(\mathbf{i})}=s_i Q_{i,i+1}^{(\mathbf{i})}
=\frac12 s_i\left((1+\mathbf{i})t_i^{(-1)}+(1-\mathbf{i})t_{i+1}^{(-1)}\right)
=\frac{1+\mathbf{i}}{2}\sigma_i+\frac{1-\mathbf{i}}{2}\sigma_i^{-1},
\end{gather*}
where $\sigma_i=s_i t_i^{(-1)}\in \mu(G)$.
Thus, $J_\mathbf{i}(R G)\subseteq R\mu(G)$, so that the automorphism $J_\mathbf{i}$ of $\mathbb C \mathbb G_n$ restricts
to an isomorphism $R G \xrightarrow{\sim} R \mu(G)$.
Theorem~\ref{thm:main2} is proved.

\subsection*{Proof of Theorem~\ref{thm:thick}}

It is convenient to prove Theorem~\ref{thm:thick} before Theorems~\ref{thm:not iso} and~\ref{thm:classification}.
We start with the following lemma.

\begin{Lemma}
\label{lem:torus}
Let~$G$ be a~thick subgroup of $G(m,1,n)={\mathcal C}^n\rtimes \mathbb S_n$ where ${\mathcal C}$ is the subgroup of
${\mathbb C}^\times$ of order $m$.
Then the group $T=G\cap {\mathcal C}^n$ is of the form $T_{{\mathcal C},{\mathcal C}'}=\{t\in {\mathcal C}^n: \det
t\in {\mathcal C}'\}$ for some subgroup ${\mathcal C}'\subset {\mathcal C}$ of $\mathbb C^\times$, and is generated by
$ \big\{t_1^{(\epsilon')}: \epsilon'\in{\mathcal C}'\big\} \cup \big\{t_i^{(\epsilon)}t_j^{(\epsilon^{-1})}: \epsilon \in {\mathcal C},1\le i,j\le n\big\}$.
\end{Lemma}

\begin{proof}
Let $\epsilon$ be a~generator of ${\mathcal C}$,  so that $t_1^{(\epsilon)}\in G(m,1,n)$.
Since~$G$ is thick, there is an element in~$G$ of the form $ts_1$ where $t\in {\mathcal C}^n$.
By the normality of~$G$, $ t_1^{(\epsilon)} (ts_1) t_1^{(\epsilon^{-1})} (ts_1)^{-1}=
t_1^{(\epsilon)}t_2^{(\epsilon^{-1})}$ belongs to~$G$, hence to $T$.
Because $\mathbb S_n$ acts on $T$ (by conjugation within~$G$), it follows that $t_i^{(\epsilon)}t_j^{(\epsilon^{-1})}\in
T$.
These elements generate the subgroup $T_0=T_{{\mathcal C},\{1\}}$ of $T$.

Every element $t'\in T$ is equal, modulo $T_0$, to an element of the form $t_1^{(\epsilon')}$ for some $\epsilon'\in
{\mathcal C}$, where $\epsilon'=\det t'$.
Denote by ${\mathcal C}'$ the group formed by all such $\epsilon'$.
Then $T\subseteq T_{{\mathcal C},{\mathcal C}'}$, and, since $\{ {t_1^{(\epsilon')}: \epsilon'\in {\mathcal C}'}\}\cup
T_0 $ generates $T_{{\mathcal C},{\mathcal C}'}$, one has $T\supseteq T_{{\mathcal C},{\mathcal C}'}$.
\end{proof}

We continue the proof of Theorem~\ref{thm:thick}.
Let~$G$ be a~thick subgroup of $G(m,1,n)={\mathcal C}^n\rtimes \mathbb S_n$ so that $G\cap {\mathcal C}^n=T_{{\mathcal
C},{\mathcal C}'}$ as in Lemma~\ref{lem:torus}.
Because~$G$ is thick,~$G$ contains an element of the form $ts_1$ where $t\in {\mathcal C}^n$.
Premultiplying $ts_1$ by an element of $T_{{\mathcal C},\{1\}}$, we conclude that $G\ni t_1^{(\epsilon)}s_1$, hence~$G$
contains $(t_1^{(\epsilon)}s_1)^2=t_1^{(\epsilon)}t_2^{(\epsilon)}$.
It follows that $t_1^{(\epsilon)}t_2^{(\epsilon)}\in T_{{\mathcal C},{\mathcal C}'}$, hence $\det
t_1^{(\epsilon)}t_2^{(\epsilon)} = \epsilon^2\in{\mathcal C}'$.
This means that
\begin{itemize}\itemsep=0pt
\item
either $\epsilon\in {\mathcal C}'$, implying $t_1^{(\epsilon)}\in G$ and $s_1\in G$;
\item
or $-\epsilon\in {\mathcal C}'$, implying $t_1^{(-\epsilon)}\in G$ and $t_1^{(-1)}s_1\in G$.
\end{itemize}
We note the following easy lemma.

\begin{Lemma}
\label{lem:obvious}
Let~$G$ be a~normal subgroup of ${\mathcal C}^n\rtimes \mathbb S_n=G(m,1,n)$ such that $G\cap {\mathcal C}^n =
T_{{\mathcal C},{\mathcal C}'}$.
\begin{enumerate}\itemsep=0pt
\item[$(a)$] If $G\ni s_1$, then~$G$ contains all elements of the form $tw$ with $t\in T_{{\mathcal C},{\mathcal C}'}$ and $w\in
\mathbb S_n$.

\item[$(b)$] If $G\ni t_1^{(-1)}s_1$, then~$G$ contains all elements $tt_1^{(\det w)}w$ with $t\in T_{{\mathcal C},{\mathcal
C}'}$ and $w\in \mathbb S_n$.
\end{enumerate}
\end{Lemma}

To continue the proof of Theorem~\ref{thm:thick}, suppose $s_1\in G$.
Then by Lemma~\ref{lem:obvious}, $G\supseteq \mathbb S_n$, Hence $G=(G\cap {\mathcal C}^n)\rtimes \mathbb S_n
=T_{{\mathcal C},{\mathcal C}'}\rtimes \mathbb S_n$.
Then $G=G(m,p,n)$ is a~complex ref\/lection group with $\frac mp=|{\mathcal C}'|$.

The only remaining case is $s_1\notin G$ but $t_1^{(-1)}s_1 \in G$.
It follows from Lemma~\ref{lem:obvious} that $t'w\in G$ (where $t'\in {\mathcal C}^n$ and $w\in \mathbb S_n$), if and
only if $t' t_1^{(\det w)}\in G$, if and only if $t' t_1^{(\det w)}\in G\cap {\mathcal C}^n = T_{{\mathcal C},{\mathcal
C}'}$, if and only if $\det \big(t' t_1^{(\det w)}\big)\in {\mathcal C}'$.
Observing that $\det \big(t' t_1^{(\det w)}\big)=\det(t'w)$, we obtain
\begin{gather*}
G = \{ t'w\in {\mathcal C}^n \rtimes \mathbb S_n: \det(t'w)\in {\mathcal C}' \} = W_{{\mathcal C},{\mathcal C}'}.
\end{gather*}
In this case, $t_1^{(-1)}s_1\in {\mathcal C}^n \rtimes \mathbb S_n$ means that $-1\in{\mathcal C}$.
That is, $m$ is even.
Theorem~\ref{thm:thick} is proved.

\subsection*{Proof of Theorem~\ref{thm:not iso}}

We will use the following notion.

\begin{Definition}
\label{def:reg}
A thick subgroup~$G$ of $G(m,1,n)$ is \textit{regular}, if for each normal abelian subgroup~$N$ of~$G$, either
$N=T_G:=G\cap (\mathbb C^\times)^n$ or $|N|<|T_G|$.
Otherwise,~$G$ is \textit{singular}.
\end{Definition}

\begin{Lemma}\label{lem:regsing}\quad
\begin{enumerate}\itemsep=0pt
\item[$(a)$] Suppose that $G_i$ is a~regular thick subgroup of $G(m_i,1,n_i)$, $i=1,2$.
If $G_1$ and $G_2$ are isomorphic (as abstract groups), then $m_1=m_2$, $n_1=n_2$ and $T_{G_1}=T_{G_2}$.

\item[$(b)$] A thick subgroup~$G$ of $G(m,1,n)$ is singular, if and only if~$G$ belongs to the following list: $G(1,1,n)$ with
$n=2,3,4$, $G(2,1,2)$, $G(2,2,2)$, $\mu(G(2,2,2))$.
\end{enumerate}
\end{Lemma}

\begin{proof}
(a) Since $G_i$ is regular and $T_{G_i}$ is the unique largest order normal abelian subgroup of~$G_i$, the restriction
of any isomorphism $f\colon G_1\to G_2$ to $T_{G_1}$ is an isomorphism $T_{G_1}\xrightarrow{\sim} T_{G_2}$, and $f$
induces an isomorphism $\overline f\colon \mathbb S_{n_1}=G_1/T_{G_1} \to \mathbb S_{n_2}=G_2/T_{G_2}$.
Hence $n_1=n_2$.
Furthermore, $m_i$ is the exponent of the group $T_{G_i}$, hence $m_1=m_2$.
Finally, $T_{G_i}=T_{{\mathcal C},{\mathcal C}_i'}$ by Lemma~\ref{lem:torus}, and $|T_{G_i}|=m_i^{n_i-1}|{\mathcal
C}_i'|$ implies ${\mathcal C}'_1={\mathcal C}'_2$.
This proves part (a).

(b) Clearly, the subgroup $T:=T_G=T_{{\mathcal C},{\mathcal C}'}$ as in Lemma~\ref{lem:torus} is a~normal abelian
subgroup of~$G$.
Let $N$ be a~normal abelian subgroup of~$G$.
We will show that $N\subseteq T$.

The map $\pi\colon G \to \mathbb S_n$ is surjective as~$G$ is thick, therefore $\pi(N)$ is a~normal abelian subgroup of
$\mathbb S_n$.
Hence if $n=1$ or $n\ge 5$, $\pi(N)$ can only be $\{1\}$ so $N\subseteq \ker \pi=T$ and~$G$ is regular.

Let $2\le n\le 4$.
Then $\mathbb S_n$ has a~unique normal abelian subgroup $K$ which is not $\{1\}$.
Assume that $N\not\subset T$.
Then $\pi(N)=K$.
One can check that $K$ acts on the indices $1,\ldots,n$ transitively, hence the centralizer of $K$ in $T$ is the set of
scalar matrices in $T$, which is the center $C(G)$ of~$G$.
Because the conjugation action of $N$ on $T$ factors through $\pi(N)$, we have $C_G(N)\cap T=C(G)$.
But $N$ is abelian, so $N\subseteq C_G(N)$.
Thus, $N\cap T \subset C(G)$.

Now let $wt\in N$ where $1\ne w\in \mathbb S_n$ and $t\in (\mathbb C^\times)^n$.
Let $i\in \{1,\ldots,n\}$ be such that $w(i)\ne i$, and let $\epsilon$ be a~generator of ${\mathcal C}$.
Then $N \ni t_i^{(\epsilon)} \cdot wt \cdot t_i^{(\epsilon^{-1})} \cdot (wt)^{-1}= t_i^{(\epsilon)}
t_{w(i)}^{(\epsilon^{-1})}$ which is a~scalar matrix only if $m=1$ or $m=n=2$.
\end{proof}

To prove the \textit{if} part of Theorem~\ref{thm:not iso}, let $m,n$ be even, $\frac mp$ be odd, and ${\mathcal C}'$,
${\mathcal C}$ be subgroups of $\mathbb C^\times$ of order $\frac mp,m$, respectively.
Assume for contradiction that there is an isomorphism ${\varphi\colon W_{{\mathcal C},{\mathcal C}'}\to G(m,p,n)}$.

In the case $G=G(2,2,2)=T_{\{\pm1\},\{1\}}\times \mathbb S_2$,~$G$ is a~Klein $4$-group which is not isomorphic to
$W_{\{\pm1\},\{1\}}$, a~cyclic group generated by $s_1\g 1$ of order $4$.

In all other cases, by Lemma~\ref{lem:torus} $T=T_{{\mathcal C},{\mathcal C}'}$ is the unique maximal normal abelian
subgroup of~$W_{{\mathcal C},{\mathcal C}'}$ and of $G=G(m,p,n)$, hence $\varphi(T)=T$ and $\varphi$ induces an
isomorphism $\overline\varphi\colon \mathbb S_n=\mu(G)/T \to \mathbb S_n=G/T$.
Let us state three easy lemmas, in which we refer to a~cycle of length $n$ in $\mathbb S_n$ as a~long cycle.

\begin{Lemma}
\label{lem:1}
Let $n$ be even and $c\in \mathbb S_n$ be a~long cycle.
Then $\det\big(\g 1 c\big)=1$ and $t_1^{(-1)} c\in W_{{\mathcal C},{\mathcal C}'}$.
\end{Lemma}

\begin{Lemma}
\label{lem:2}
For $t\in (\mathbb C^\times)^n$ and a~long cycle $c$, $(tc)^n=z^{(\det t)}$, where $z^{(\epsilon)}$ denotes
$t_1^{(\epsilon)}t_2^{(\epsilon)}\cdots t_n^{(\epsilon)}$.
\end{Lemma}

\begin{Lemma}\label{lem:3}
The image of a~long cycle in $\mathbb S_n$ under any automorphism of $\mathbb S_n$ is a~long cycle.
\end{Lemma}

Let $c$ be a~long cycle in $\mathbb S_n$.
By Lemma~\ref{lem:1}, $\overline c=t_1^{(-1)}c\in\mu(G)$.
By Lemma~\ref{lem:2}, $(\overline{c})^n=z^{(-1)}$.
Hence $(\overline{c})^{n|{\mathcal C}'|}=\big(z^{(-1)}\big)^{|{\mathcal C}'|}= z^{(-1)}$ as $|{\mathcal C}'|=\frac mp$ is odd.

Now consider $\varphi(\overline c)\in G$.
It is of the form $t\overline{\varphi}(c)$, where $t$ is some element of ${\mathcal C}^n$ such that, by def\/inition of
$G=G(m,p,n)$, $\det t\in {\mathcal C}'$.
By Lemma~\ref{lem:3}, $\overline{\varphi}(c)$ is a~long cycle in $\mathbb S_n$.
Therefore, by Lemma~\ref{lem:2},
\begin{gather*}
(\varphi(\overline c))^{n|{\mathcal C}'|} = \big(z^{(\det t)}\big)^{|{\mathcal C}'|} = z^{((\det t)^{|{\mathcal C}'|})} =1.
\end{gather*}
This is a~contradiction, because the image of $z^{(-1)}\ne 1$ under an isomorphism $\varphi$ cannot be $1$.
The \textit{if} part of Theorem~\ref{thm:not iso} is proved.

To establish the \textit{only if} part of Theorem~\ref{thm:not iso}, observe that if $\frac mp$ is even, then the groups
$G=G(m,p,n)$ and $\mu(G)=W_{{\mathcal C},{\mathcal C}'}$ simply coincide as subgroups of $\mathbb G_n$.
Indeed, in this case $\det w\in \{\pm1\}\subseteq {\mathcal C}'$, hence the conditions $\det t\in {\mathcal C}'$ and
$\det tw\in {\mathcal C}'$ are equivalent.

If $m$ is even, $\frac mp$ is odd and $n$ is odd, both $G=G(m,p,n)$ and $\mu(G)$ are normal subgroups of $\widetilde
G=G(m,\frac p2,n)$ which do not contain the central subgroup $Z=\big\{1,z^{(-1)}\big\}$ of $\widetilde G$.
Hence $\widetilde G = G\times Z = \mu(G)\times Z$, so that $G\cong\mu(G)\cong \widetilde G / Z$.
Theorem~\ref{thm:not iso} is proved.

\subsection*{Proof of Theorem~\ref{thm:classification}}

Let $G\subseteq G(m,1,n)$, $G'\subseteq G(m',1,n')$ be thick subgroups.

Assume that $G\ne G'$ are regular in the sense of Def\/inition~\ref{def:reg}.
Then by Lemma~\ref{lem:regsing}(a), $G\cong G'$ implies $n=n'$, $m=m'$ and $T_G=T_{G'}$, so we are done by
Theorem~\ref{thm:thick}.

By inspection of the list in Lemma~\ref{lem:regsing}(b), no two singular subgroups are isomorphic to each other.
Also by inspection, one shows that if $n=n'$,~$G$ is singular and~$G'$ is regular, then~$G$ and~$G'$ are never
isomorphic.

The only remaining case is when $n\ne n'$,~$G$ is singular and $G'$ is regular.
If $G\cong G'$, then $n,n'\le 4$, because $|\mathbb S_5|>|G|$.
Then, an easy analysis based on the cardinalities of singular groups shows that the only possible isomorphism is the one
given in part~(b).
Theorem~\ref{thm:classification} is proved.

\section{Appendix}

The aim of this section is to prove the following important result about group ring actions.

\begin{Theorem}
\label{thm:Premet}

Let $A$ be an integral domain and~$G$ be a~group of ring automorphisms of $A$.
Then:
\begin{enumerate}\itemsep=0pt
\item[$(a)$] the natural map $\rho\colon AG \to \End_\mathbb Z(A)$ given by $\big(\rho\big(\sum\limits_i a_i g_i\big)\big)(a)=\sum\limits_i a_i g_i(a)$ is injective;
\item[$(b)$] with respect to the natural ring structure on $\End_\mathbb Z(A)$ and the semidirect product structure $A\rtimes
\mathbb Z G$ on $AG$, the map $\rho$ is a~ring homomorphism.
\end{enumerate}
\end{Theorem}

The following is immediate from the theorem.

\begin{Corollary}
\label{cor:premet-dedekind}
In the notation of the theorem, let $R$ be a~subring of $A^G$.
Then the restriction of $\rho$ to $RG$ is an injective ring homomorphism
\begin{gather*}
RG\hookrightarrow \End_R A,
\end{gather*}
where $RG$ is the ordinary group ring of~$G$.
\end{Corollary}

\begin{proof}[Proof of Theorem~\ref{thm:Premet}]
We need the following generalization of the celebrated Dedekind's lemma.
\begin{Lemma}
\label{lem:dedekind}
Let $A$ be an integral domain and $B$ be a~multiplicative monoid.
Let $\mathcal G$ be a~set of monoid homomorphisms $B\to A$.
Then the natural $A$-linear map
\begin{gather*}
\rho\colon \ A\mathcal G \to {\rm Fun}(B,A),
\qquad
\left(\rho\left(\sum\limits_g a_g g\right)\right)(b) = \sum\limits_g a_g g(b)
\end{gather*}
is injective.
$($Here $A\mathcal G=\oplus_{g\in \mathcal G} Ag$ is the free $A$-module generated by $\mathcal G.)$
\end{Lemma}
\begin{proof}
Assume for contradiction that $\rho$ is not injective.
Let $\mathcal G _0$ be a~minimal f\/inite subset of $\mathcal G $ such that there exists a~non-zero element
$\sum\limits_{g\in \mathcal G _0}a_g g$ in the kernel of $\rho$.
That is,
\begin{gather*}
\sum\limits_{g\in \mathcal G _0}a_g g(b)=0
\qquad
\text{for all} \ \ b\in B.
\end{gather*}
Clearly, all $a_g$ are non-zero and $|\mathcal G _0|>1$ because $A$ has no zero divisors.
In particular, for each $b,b'\in B$,
\begin{gather*}
\sum\limits_{g\in \mathcal G _0} a_g g(b'b)=\sum\limits_{g\in \mathcal G _0} a_g g(b')g(b)=0,
\end{gather*}
because all $g\in \mathcal G _0$ are homomorphisms from $B$ to $A$.
Furthermore, f\/ix $h\in \mathcal G _0$.
Combining the above identities, we obtain
\begin{gather*}
\sum\limits_{g\in \mathcal G _0\setminus\{h\}} a_g g(b') (g(b) - h(b)) = 0
\qquad
\text{for all}\quad
b,b'\in B.
\end{gather*}
That is, for each $b\in B$, the element
\begin{gather*}
k_b:=\sum\limits_{g\in \mathcal G _0\setminus\{h\}} a_g (g(b) - h(b)) g
\end{gather*}
belongs to $\ker \rho$.
The minimality of $\mathcal G _0$ implies that $k_b=0$ for all $b\in B$, which is equivalent to $g(b)=h(b)$ for all
$b\in B$ and all $g\in \mathcal G _0$.
That is, $|\mathcal G _0|=1$, a~contradiction.
\end{proof}

We now use Lemma~\ref{lem:dedekind} for $A$ as in Theorem~\ref{thm:Premet} with $B=A\setminus\{0\}$ and $\mathcal G=G$
viewed as homomorphisms of monoids $A\setminus\{0\} \to A$.
Since $\End_\mathbb Z A$ is naturally a~subset of ${\rm Fun}(A\setminus\{0\},A)$ and $\rho(AG)\in \End_\mathbb Z
A\subset{\rm Fun}(A\setminus\{0\},A)$, part (a) of Theorem~\ref{thm:Premet} is proved.

To prove (b), note that $\End_\mathbb Z A$ is naturally a~ring, with multiplication being composition of maps.
The semidirect product multiplication on $A G$ is given by the formula $(ag)(a'g')=ag(a')\cdot gg'$ for all $a,a'\in A$,
$g,g'\in G$.
Then
\begin{gather*}
\rho((ag)(a'g'))(b)=a\cdot g(a'\cdot g'(b))=a\cdot g(a')\cdot g(g'(b)) = ag(a') (gg')(b)
\end{gather*}
for all $a,a',b\in A$, $g,g'\in G$.
Part (b) of the theorem is proved.
\end{proof}

\subsection*{Acknowledgments}

We thank Ken Brown for bringing the paper~\cite{KKZ} to our attention, and Alexander
Premet for stimulating discussions.
The present paper was started when both authors were research members of the Mathematical Sciences Research Institute.
We thank the Institute and the organizers of the Noncommutative Algebraic Geometry and Representation Theory program for
creating an atmosphere conducive for research.
We acknowledge partial support of the LMS Research in Pairs grant ref.~41224.
The second named author was partially supported by the NSF grant DMS-1101507.

\pdfbookmark[1]{References}{ref}
\LastPageEnding

\end{document}